\documentclass[12pt]{article}
\usepackage{amssymb}
\begin{document}
\begin{center}
{\Large Some Characterizations of VNL Rings}\\[3mm] {\small
Harpreet K.Grover\footnote{The research of the first author is
supported by CSIR, India and will form part of her Ph.D Thesis} and
Dinesh Khurana\footnote{Corresponding author}}\\[3mm]
{\footnotesize Department of Mathematics, Panjab University,
Chandigarh-160014, India\\ \em
harpreetgr@gmail.com, dkhurana@pu.ac.in}
\end{center}

\vspace{5mm} \noindent{\bf \footnotesize Abstract}\\[3mm]
{\footnotesize A ring R is said to be VNL if for any $a\in R$,
either $a$ or $1-a$ is (von Neumann) regular. The class of VNL
rings lies properly between the exchange rings and (von Neumann)
regular rings. We characterize abelian VNL rings. We also
characterize and classify arbitrary VNL rings without infinite set
of orthogonal idempotents; and also the VNL rings having primitive
idempotent $e$ such that $eRe$ is not a division ring. We prove that a semiperfect
ring $R$ is VNL if and only if
for any right uni-modular row $(a_1, a_2)\in R^2$, one of the $a_i$ is regular in $R$.
  Formal
triangular matrix rings that are VNL, are also characterized.
As a corollary it is shown that an upper
triangular matrix ring $T_n(R)$ is VNL if and only if $n=2$ or $3$
and $R$ is a division ring. }

\vspace{7mm}
\noindent{\bf 1. Introduction}

\vspace{3mm}As a common generalization of local and (von Neumann) regular rings,
 Contessa in [5] called a ring $R$ VNL (von Neumann local) if for
each $a\in R$, either $a$ or $1-a$ is (von Neumann) regular. As
every regular element $a$ of a ring $R$ is an exchange element (in
the sense that there exists an idempotent $e\in aR$ such that $1-e
\in (1-a)R$), VNL rings are exchange rings. But if $R$ is a local
ring with nonzero $J(R)$,
 then $R\times R$, which is an exchange ring, is not a VNL ring. For instance, if $a=(x, 1-x)$, where $x$ is a nonzero
 element in $J(R)$, then neither $a$ nor $1-a$ is regular. Although VNL rings have been studied in some detail (see [3], [4], [8] and [9]), their structure is not known even in commutative case. For instance, Osba, Henriksen and Alkam in ([8], page 2641) remark:

\vspace{2mm}
  {\em We are unable to characterize (commutative) VNL-rings abstractly in the sense of relating them to more familiar classes of rings...}

\vspace{2mm}
   The present paper is an effort towards this direction. We characterize abelian\footnote{A ring is called abelian if its all idempotents are central.} VNL rings. It is shown in Section 3 below that abelian VNL rings are precisely those exchange rings $R$ in which, for any idempotent $e$,  one of the two corner rings $eRe$ and $(1-e)R(1-e)$ is regular. Let $M(R)$ denote the maximal regular ideal of a ring $R$ as defined by Brown and McCoy in [1]. In ([3], Lemma 2.7) Chen and Tong showed that if $R$ is an abelian VNL ring, then $R/M(R)$ is a local ring. We show that abelian VNL rings are precisely those rings $R$ for which $R/M(R)$ is a local ring.  But this characterization of abelian VNL rings is not valid for arbitrary rings (see Example 3.3 below). In Section 4, we characterize arbitrary VNL rings which do not have infinite set of orthogonal idempotents. As an exchange ring without infinite set of orthogonal idempotents is semiperfect (see [2]),  this gives us characterization of semiperfect VNL rings. We prove that a semiperfect ring $R$ is VNL if and only if for any right uni-modular row $(a_1, a_2)\in R^2$, one of the $a_i$ is regular in $R$. We also characterize VNL  rings $R$ with a primitive idempotent $e$ such that $eRe$ is not a division ring
   or equivalently $J(eRe)\ne 0$ (if $e$ is a primitive idempotent in an exchange ring $R$, then $eRe$ is a local ring).

   In Section 2, we give some examples of VNL rings and prove some basic properties of VNL rings. For
   instance, it is proved that if $e$ is an idempotent in a VNL ring $R$, then either $eRe$ or $(1-e)R(1-e)$ is a regular ring. We also show that this property does not characterize VNL rings. We also understand the regular elements of formal triangular matrix ring $\left(\begin{array}{cc} R & M\\ 0 & S \end{array}\right)$ and with the help of this, we characterize formal triangular matrix rings, that are VNL. As a corollary, we prove that the upper
triangular matrix ring $T_n(R)$ is VNL if and only if $n=2$ or $3$
and $R$ is a division ring.

In [6] Nicholson defined a ring $R$ to be NJ if every element of $R\verb=\=J(R)$ is regular. Clearly, an NJ ring $R$ is VNL and Nicholson proved that in an NJ-ring, $eRe$ is regular for every proper idempotent $e$ of $R$. In Section 5, we prove that if $R$ is a ring without a nontrivial central idempotent and $J(R)\ne 0$, then $R$ is NJ if and only if $R$ is VNL and $J(eRe)=0$ for every proper idempotent $e$ of $R$.\\[6mm]
 \noindent {\bf 2. Examples and Basic Properties of VNL rings}

 \vspace{3mm}
 The trivial examples of VNL rings, of course, are
regular and local rings. Here we give some non-trivial examples of
VNL rings.\\[4mm]\noindent {\bf Examples 2.1}\\[3mm]
 \noindent (1) Let $R =\{(q_1,q_2,...,q_n,z,z,z,...):n\geq 1, q_i
\in\mathbb{Q}~ and ~  z\in\mathbb{Z}_2\}$ where $\mathbb{Z}_2$
denotes the localization of $\mathbb{Z}$ at the prime ideal $(2)$.
Then $R$ is a VNL ring. An element $(q_1, q_2,...,q_n, z,z,z,
\dots)$ is regular precisely when $z$ is a unit in $\mathbb{Z}_2$.
It is easy to see that  every non-zero ideal of $R$ contains a
non-zero idempotent implying that $J(R)=0$. Thus VNL rings may not
be semiregular.\\[2mm]\noindent (2) Nicholson in [6] studied the
rings $R$ with the property that every element outside $J(R)$ is
regular. He called these rings NJ-rings. Clearly, every NJ ring is
VNL. He characterized NJ rings by showing that the only NJ rings
besides regular and local rings are of the form
$\left(\begin{array}{cc} D_1 & X\\ Y & D_2 \end{array}\right)$,
where $D_1$ and $D_2$ are division rings, $X$ and $Y$ respectively
are $(D_1, D_2)$ and $(D_2, D_1)$ bimodules with $XY=0=YX$. If $D$
is a division ring, then from Nicholson's characterization of NJ
rings, it is clear that the upper triangular matrix ring $T_2(D)$
is an NJ ring and hence VNL.
\\[2mm] \noindent (3) It was observed in [8] that the ring
$\mathbb{Z}_n$ of integers mod $n$ is VNL if and only if $(pq)^2$
does not divide $n$ where $p$ and $q$ are distinct primes. This
is clear from the fact that if $R\times S$ is a VNL ring, then either $R$ or
$S$ is regular. \\[2mm] \noindent (4) For a commutative ring $R$, the formal
power series ring $R[[x]]$ is VNL if and only if $R$ is local (see
[8]).
\\[2mm] \noindent (5) If $R$ is a regular ring, $L$ is a local
ring and $_R M_L$ is a bimodule then $\left(\begin{array}{cc} R &
M\\ 0 & L
\end{array}\right)$ is VNL. In fact, we show that every element of
the type $\left(\begin{array}{cc} r & m\\ 0 & l
\end{array}\right)$, where $l$ is a unit in $L$, is regular. Since
$r$ is regular so there exists an $s\in R$ such that $rsr=r$, then as $l$ is a
unit in $L$, so
 $$\left(\begin{array}{cc} r & m\\ 0 & l
\end{array}\right) \left(\begin{array}{cc} s & -sml^{-1}\\ 0 &
l^{-1}
\end{array}\right) \left(\begin{array}{cc} r & m\\ 0 & l
\end{array}\right) = \left(\begin{array}{cc} r & m\\ 0 & l
\end{array}\right).$$

\vspace{4mm}
We now prove some basic results about VNL rings. \\[4mm]
 {\bf Proposition 2.2.} {\em Let $R$ be a VNL ring then center of
$R$ is also a VNL ring.}\\[3mm] {\bf Proof.} Let $x\in Z(R)$ be
regular then there exists $y\in R$ such that $x=xyx$.Then it is easy to
see that $z=yxy\in Z(R)$ and $x=xzx$ i.e $x$ is regular in $Z(R)$.
So center of $R$ is also VNL. $\Box$\\[2mm] The following corollary
is immediate from above result:\\[3mm] {\bf Corollary 2.3.} {\em
Let $R$ be a VNL ring, then it is indecomposable as a ring if and
only if its center is a local ring.}\\[3mm]\indent If $R = S\times
T$ is VNL,
 then it
is clear that either $S$ or $T$ is a regular ring, because if $s$ in $S$
and $t$ in $T$ are non-regular elements
then neither $r=(s, 1-t)$ nor $1-r = (1-s, t)$ is regular. Thus if
$e$ is a central idempotent in a VNL ring $R$, then either $eRe$
or $(1-e)R(1-e)$ is regular. Interestingly, this also
holds for non-central idempotents of VNL rings as shown below.\\[4mm]
 {\bf Lemma 2.4.} {\em If $R$ is a VNL ring then for
every idempotent $e$ of $R$, either $eRe$ or $(1-e)R(1-e)$ is a
regular ring.}\\[3mm] {\bf Proof.} Let $R$ be a VNL ring and $e\in
R$ be an idempotent. Then $$R\cong \left(\begin{array}{cc} eRe &
eR(1-e)\\ (1-e)Re & (1-e)R(1-e)\end{array}\right).$$ If $ x\in
eRe$ and $y\in~(1-e)R(1-e)$ are two non-regular elements, then
both $a=$ $\left(\begin{array}{cc} x & 0\\ 0 & 1-y
\end{array}\right)$ and $1-a=$ $\left(\begin{array}{cc} 1-x & 0\\
0 & y
\end{array}\right)$ are also non-regular. $\Box$

\vspace{3mm} The following example shows that the necessary
condition of Lemma 2.4 above does not characterize VNL
rings.\\[4mm]{\bf Example 2.5.} Let  $R=$
$\{(q_1,q_2,....,q_n,z,z,z,....): n\geq 1, q_i\in\mathbb{Q},
z\in\mathbb{Z} \}$. It is clear that for idempotent $e$ of $R$,
either $eRe$ or $(1-e)R(1-e)$ is regular. But $R$ is not a VNL
ring as $\mathbb{Z}$ is a homomorphic image of $R$, which is not
VNL.

\vspace{3mm} The following result also clearly follows from Lemma
2.4.\\[3mm] {\bf Corollary 2.6.} {\em For a ring $R$, the matrix
ring $M_n(R)$, $n>1$, is VNL if and only if $R$ is regular.}

\vspace{3mm} In [8], Osba, Henriksen and Alkam defined a
commutative ring $R$  to be SVNL if $\sum_{i=1}^{n}a_iR=R$ implies
that one of the $a_i$'s is regular, and asked if every commutative
VNL ring
  $R$ is SVNL. This question was answered by Chen and Tong in [3], where they, in fact, proved that
  whenever $\sum_{i=1}^{n}a_iR=R$ in an abelian VNL ring, then one of the $a_i$'s is regular. We give below a different
  proof of their result.\\[3mm]{\bf Corollary 2.7} (Chen and Tong [3], Theorem 2.8). {\em  Let $R$ be an
abelian VNL ring. If $\sum_{i=1}^{n}a_iR=R$, then one of the $a_i$'s is regular.}\\[2mm]{\bf Proof.}
 As $R$ is an exchange ring and $\sum_{i=1}^{n}a_iR=R$,
 there exist an orthogonal set $\{e_1, \dots, e_n\}$ of idempotents such that
 $e_i \in a_iR$ and $e_1+\dots+ e_n = 1$ (see [7], Proposition 1.11).
Now for each $i$, $$e_i\in a_i R \Longrightarrow a_i R + (1-e_i) R
= R \Longrightarrow e_i a_i R = e_i R. $$  Thus $e_i a_i$ is
regular for every $i$.  By Lemma 2.4, either $e_iR=e_i R
e_i$  or $(1-e_i)R=(1-~e_i)~R~(1-~e_i)$ is regular for each $i$. If $e_iRe_i$
is regular for each $i$, then $R$, being a direct product of
$e_iRe_i$, is regular. Also if $(1-e_i)R(1-e_i)$ is regular for
some $i$, then, as $e_ia_i$ is already regular,
$a_i=e_ia_i+(1-e_i)a_i$ is regular. $\Box$

\vspace{3mm}
We now characterize the regular elements of  formal triangular
matrix rings. The characterization turned out to be very useful in the investigation of VNL rings.\\[3mm]{\bf Proposition 2.8.} {\em Let  $_RM_S$ be a bimodule. An element
$\left(\begin{array}{cc} a & m\\ 0 & b
\end{array}\right)$ of the formal triangular matrix ring $T=$
$\left(\begin{array}{cc} R & M\\ 0 & S \end{array}\right)$ is
regular if and only if there exist idempotents $e\in R$ and $f\in S$ such that $aR=eR$, $Sb=Sf$ and $(1-e)m(1-f)=0$.}\\[3mm]
{\bf Proof.} If $\left(\begin{array}{cc} a & m\\ 0 & b
\end{array}\right)$ is regular in $T$, then  for some
$\left(\begin{array}{cc} x & y\\ 0 & z \end{array}\right)$ in $T$
$$\left(\begin{array}{cc} a & m \\ 0 & b
\end{array}\right) \left(\begin{array}{cc} x & y\\ 0 & z
\end{array}\right) \left(\begin{array}{cc} a &  m\\ 0 & b
\end{array}\right) = \left(\begin{array}{cc} a & m\\ 0 & b
\end{array}\right).$$ So
$ axa=a,\;\;  bzb=b,\;\; axm+ayb+mzb=m.
$ If we take $e=ax$ and $f=zb$, then $$aR=eR,\;\; Sb=Sf \mbox{ and }
(1-e)m(1-f)=0.$$
Conversely, let $aR=eR$, $Sb=Sf$ and
$(1-e)m(1-f)=0$ for some idempotents $e\in R$ and $f\in S$. Then $m=em+mf-emf$, $ar=e$ and $sb=f$ for some $r\in R$ and $s\in S$. Thus
$$\left(\begin{array}{cc} a & m\\ 0 & b \end{array}\right)
\left(\begin{array}{cc} r & -rms\\ 0 & s \end{array}\right)
\left(\begin{array}{cc} a & m\\ 0 & b
\end{array}\right)=\left(\begin{array}{cc} a & m\\ 0 & b
\end{array}\right).\;\;\Box$$

\vspace{4mm}
For the characterization of
 formal triangular matrix rings that are VNL,  we need
the following\\[3mm] {\bf Definition 2.9.} We call a module $_RM$
{\em partial} if for any idempotent $e\in R$, either $eM=0$ or
$(1-e)M=0$.

\vspace{4mm} Here are some examples of partial modules.\\[3mm]
{\bf Examples 2.10.} (1) Any module over a ring $R$ with only
trivial idempotents is partial. In particular, every vector space
is a partial module.\\[2mm] (2) Any simple module over a
commutative ring $R$ is partial.\\[2mm] (3) If $S$ is a ring with
only trivial idempotents, $R = S\times S\times\dots \times S$ and
$M = S\times 0\times\dots \times 0$, then clearly $_R M$ is a
partial module.

\vspace{3mm} It follows from the following result that no non-zero
module over a proper matrix ring is partial.\\[3mm] {\bf
Proposition 2.11.} {\em For any ring $R$, let $S=M_n(R)$ with
$n\geq 2$, and $_S M$ be a non-zero module. Then for $0\neq m\in
M$, there exists an idempotent $e$ in $S$ such that $em\neq 0$ and
$(1-e)m\neq 0$. In particular, no non-zero module over $S$ is
partial.}\\[2mm]{\bf Proof.}  In $S$ $$1=E_{11}+E_{22}+....+E_{nn}
\Rightarrow m = E_{11}m+E_{22}m+....+E_{nn}m.$$  Since $m$ is
non-zero,
 $E_{ii}m\neq 0$ for some $i$. If for some $j\neq i$,
$E_{jj}m\neq 0$, then either $(1-E_{ii})m\neq 0$ or
$(1-E_{jj})m\neq 0$, because otherwise $m = E_{ii}m = E_{jj}m$ implying that $m = 0$.
 Now suppose that there is only one $i$ such that
$E_{ii}m\neq 0$. Pick any
$j\neq i$ and consider $e = E_{jj}+E_{ji}$. Then $e$ is an idempotent in $S$. Now
$em = E_{jj}m+E_{ji}m = 0+E_{ji}m = E_{ji}m$. If $E_{ji}m = 0$,
then $ E_{ij}E_{ji}m = 0$ implying that $E_{ii}m = 0$. Hence $em\neq 0$. Also $$(1-e)m = m-E_{jj}m-E_{ji}m =
m-E_{ji}m.$$ If $m-E_{ji}m = 0 $, then $E_{ji}m-E_{ji}^2m = E_{ji}m = 0$. But this, as seen above, implies that $m=0$. $ \Box$

\vspace{3mm} In the following result we characterize the formal
triangular matrix rings that are VNL.\\[4mm] {\bf Theorem 2.12.}
{\em Let $_RM_S$ be a bimodule. Then the ring $T =$
$\left(\begin{array}{cc} R & M\\ 0 & S \end{array}\right)$ is VNL
if and only if \\[2mm] (1) One of $R$ and $S$ is regular and the
other is VNL.\\[2mm] (2) Either $_R M$ or $M_S$ is a partial
module.\\[2mm] (3) For any non-regular $r\in R$, $(1-r)M=M$ and
for any non-regular element
 $s\in S$, $M(1-s)=M$.}\\[2mm] {\bf Proof.}
Suppose that $T$ is VNL. Then, by Lemma 2.4 and the fact that every factor ring of
 a VNL ring is VNL, one of the $R$ and $S$ is regular
and the other is VNL.  Now suppose that $_R M$
is not a partial module. So there exist an idempotent $e\in R$
such that $$ eM\neq 0 \mbox { and }
  (1-e)M\neq 0.\hspace{10mm} (A)$$ Let $f$ be an idempotent in $S$, then if we take idempotent $E= \left(\begin{array}{cc} e & 0\\ 0 & f \end{array}\right) \in T$,
 by Lemma 2.4, either $ETE$ or $(1-E)T(1-E)$ is regular. So
$$ \mbox{ either } eMf = 0 \mbox{ or } (1-e)M(1-f) =
0.\hspace{10mm} (B)$$ Similarly, if we take the idempotent
$\left(\begin{array}{cc} 1-e & 0\\ 0 & f \end{array}\right)\in T$,
we get $$\mbox{ either } (1-e)Mf = 0 \mbox{ or } eM(1-f) =
0.\hspace{10mm} (C)$$ From (A), (B) and (C), it is clear that
either $Mf = 0 $ or $M(1-f) = 0$. Thus\ $ M_S$ is a partial
module. Now suppose that $r$ is a non-regular element in $R$. Then
as  $\left(\begin{array}{cc} r & m\\ 0 & 1 \end{array}\right)$ is
not regular for any $m\in M$ , the element
$\left(\begin{array}{cc} 1-r & m\\ 0 & 0
\end{array}\right)$ is regular for every $m\in M$. Now if $(1-r)R
= eR$, then by Proposition 2.8, $(1-e)M = 0$. So $M= eM =(1-r)M$.
Similarly if $s$ is a non-regular element of $S$, then $M(1-s) =
M$.

 For converse,  we may assume without loss of generality that $R$ is VNL and $S$ is regular. Let $x=$ $\left(\begin{array}{cc} r & m \\ 0 & s \end{array}\right)$ $\in T$. If $r$ is regular in $R$ but $1-r$ is not regular, then we show that $x$ is regular in $T$. As $1-r$ is not regular, by condition (3), $rM = M$. Thus if $rR = eR$ for some idempotent $e\in R$, then $eM = M$ and so $(1-e)M = 0$. This, by Proposition 2.8, implies that $x$ is regular. Now suppose that both $r$ and $1-r$ are regular in $R$. Suppose $_R M$ is a partial module. Now if $rR = eR$ for some idempotent $e\in R$, then either $eM = 0$ or $(1-e)M = 0$. If $(1-e)M = 0$, then, by Proposition 2.8,  $x$ is regular. Now if $eM = 0$, then $rM = 0$. So $(1-r)M = M$ and if $(1-r)R = fR$ for some idempotent $f$ of $R$, then $fM = M$ implying that $(1-f)M = 0$. Thus, by Proposition 2.8, $1-x$ is regular. So $T$ is VNL. Lastly if $M_S$ is partial, we can similarly prove that $T$ is VNL. $\Box$

\vspace{5mm}
 We now give various applications of Theorem 2.12. In [4], it was proved that
 if $D$ is a division ring, then $T_2(D)$ and $T_3(D)$ are VNL. The following characterization shows that
 these are the only upper triangular matrix rings that are VNL.\\[5mm]
 {\bf Corollary  2.13.} {\em The upper triangular matrix ring $T_n (R)$ is VNL
  if and only if $n=2$ or $3$ and $R$ is a division
ring.}\\[3mm]{\bf Proof.}  Let $n\geq 4$ and $e = E_{11}+E_{22}$.
Then $e$ is an idempotent in $T_n (R)$ and $eT_n(R)e\cong T_2(R)$,
 $(1-e)T_n(R)(1-e)\cong T_{n-2}(R)$ are both not regular. So by
Lemma 2.4, $T_n(R)$ is not regular if $n\geq 4$.  Also it is clear
by Proposition 2.8 that any element outside the Jacobson radical
of $T_2(D)$ is regular implying that $T_2(D)$ is VNL. Now
$T_3(D)=\left( \begin{array}{cc}T_2(D) & M\\0 & D \end{array}
\right)$. Clearly, $T_3(D)$ satisfies the conditions (1) and (2)
of Theorem 2.12. Also
 if $r\in T_2(D)$ is non-regular, then $1-r$ is a unit and so $(1-r)M=M$, implying that $T_3(D)$ also satisfies the condition (3) of Theorem 2.12 and is thus VNL.

\vspace{1mm}
 Now suppose that $T_2(R)$ is VNL.  By Lemma 2.4, $R$ is regular. If $e$ is any non-trivial idempotent
 in $R$, then neither $_RR$ nor $R_R$ is partial. So, by Theorem 2.12, $T_(R)$ is not VNL. Thus $R$ is a regular ring
 without non-trivial idempotents and is thus a division ring. Lastly if $T_3(R)$ is VNL, then so is $T_2(R)$ being
  a homomorphic image of $T_3(R)$, implying again that $R$ is a division ring. $\Box$\\[5mm]
 {\bf Corollary 2.14.} {\em (1) If $R$ is a regular ring and $S$ is a local ring, then for any bimodule $_RM_S$, the ring $\left( \begin{array}{cc} R & M\\0&S \end{array} \right)$ is a VNL ring.\\[2mm]
 (2) For any division ring $D$, the ring $\left( \begin{array}{cc} M_2(D) & M\\0&T_2(D) \end{array} \right)$, where
 $M=\left( \begin{array}{cc} 0 & D\\0&D \end{array} \right)$, is VNL.}\\[5mm]
 {\bf Proof.} The part (1) is immediate from Theorem 2.12. Also, as $M_{T_2(D)}$ is partial and for any non-regular
 $x\in T_2(D)$, $1-x$ is a unit, $\left( \begin{array}{cc} M_2(D) & M\\0&T_2(D) \end{array} \right)$ satisfies
 all the three conditions of Theorem 2.12 and is thus VNL. $\Box$\\[5mm]
  {\bf Corollary 2.15.} {\em Let $R$
and $S$ be simple artinian rings and $_R M_S$ be a bimodule. Then
the ring $T=$ $\left(\begin{array}{cc} R & M\\ 0 & S
\end{array}\right)$ is VNL if and only if either $M=0$ or one of $R$ and $S$
is a division ring.}\\[2mm]{\bf Proof.} The `if' part follows from
Corollary 2.14(1). Conversely, suppose $T$ is VNL. If $M$ is
non-zero and neither $R$ nor $S$ is a division ring, then  $_R M$
and $M_S$ are not partial by Proposition 2.11. This, in view of
Theorem 2.12 implies that $T$ is not VNL. $\Box$

\vspace{3mm} The above proof, in fact, shows
if $_{M_{n} (R)}M_{M_m (S)}$ is a bimodule, then the ring
$\left(\begin{array}{cc} M_n(R) & M\\ 0 & M_m(S)
\end{array}\right)$ is VNL only if either $M=0$ or either $n=1$ or $m=1$.\\[4mm]
 {\bf Corollary 2.16.} {\em Let $R$ be a
commutative ring and $I$ be non-zero ideal of $R$. Then $T=$
$\left(\begin{array}{cc} R & I\\ 0 & R \end{array}\right)$ is VNL
if and only if $R=F\times S$ where $F$ is a field and $S$ is a
regular ring and $I=F\times 0$.}\\[3mm]{\bf Proof.} The
sufficiency follows easily from Theorem 2.12. Conversely suppose
$T=$ $\left(\begin{array}{cc} R & I\\ 0 & R
\end{array}\right)$ is VNL. Then by Theorem 2.12, $R$ is
regular and for any idempotent $e$ of $R$, either $eI = 0$ or
$(1-e)I = 0$.  Now for any $0\neq a\in I$, $aR=eR$ for some
idempotent $0\neq e$ in $R$. As $eI\neq 0$, $(1-e)I=0$ and so
$I\subseteq eR=aR\subseteq I$. Thus $I=eR=aR$ for every nonzero
$a\in I$. This also shows that for any $a\in I$, $aI=I$, implying
that $I=eR$ is a simple, commutative, regular ring. Thus $I$ is a
field and $R=I \times (1-e)R$. $\Box$\\[4mm]\indent We now show
that the previous result also holds for non-commutative rings that
do not have infinite set of orthogonal idempotents.\\[4mm]{\bf
Corollary 2.17.} {\em Let $R$ be a ring which does not have
infinite set of orthogonal idempotents and let $I$ be a non-zero
two sided ideal of $R$. Then $\left(\begin{array}{cc} R & I\\ 0 &
R
\end{array}\right)$ is VNL if and only if $R=D\times S$ where $D$ is a
division ring, $S$ is a semisimple ring and $I=D\times 0$.}\\[3mm]
{\bf Proof.} Sufficiency is clear from Theorem 2.12. If $T=$
$\left(\begin{array}{cc} R & I\\ 0 & R \end{array}\right)$ is VNL
then $R$ is regular by Lemma 2.4. As $R$ does not have infinite set of
orthogonal idempotents, $R$ is semisimple. Thus $R =
M_{n_1}(D_1)\times M_{n_2}(D_2)\times \dots \times M_{n_r}(D_r)$,
where $D_i$'s are division rings. So $I = I_1\times I_2\times
\dots\times I_r$ where each $I_i = 0$ or $M_{n_i}(D_i)$. For any
central idempotent $e\in R$, by Theorem 2.12, either $eI=0$ or
$(1-e)I=0$. So it is clear that exactly one $I_i$ is non-zero and
we may assume $I_1$ is non-zero.  So $T=\left( \begin{array}{cc}
M_{n_1}(D_1)\times S & M_{n_1}(D_1)\times 0\\0& M_{n_1}(D_1)\times
S
\end{array} \right)$. Now
$E=\left(\begin{array}{cc} (1,0) & (0,0)\\ 0 & (1,0)
\end{array}\right)$ is a central idempotent in $T$ and so  $ET =
\left(\begin{array}{cc} M_{n_1}(D_1) & M_{n_1}(D_1)\\ 0 &
M_{n_1}(D_1) \end{array}\right)$ is VNL. So by Corollary 2.13,
$n_1=1$. $\Box$

\vspace{3mm} We will need the following result, which was proved
by Chen and Ying in [4]. As the paper is yet to appear, we give
their proof below.\\[4mm] {\bf Lemma 2.18.} {\em If $R$ is VNL
then so is $eRe$ for any idempotent $e$ in $R$.}\\[3mm] {\bf
Proof.} Let $a\in eRe$. Suppose $a$ is not regular in $R$. Then
$1-a$ is regular. Suppose $1-a=(1-a)b(1-a)$, for some $b\in R$.
Then $e-a=e(1-a)e=e(1-a)b(1-a)e=(e-a)ebe(e-a)$. $\Box$ \\[8mm]
{\bf 3. Characterizations of abelian VNL rings}

\vspace{5mm}
In this section we characterize abelian VNL rings.
\\[4mm]{\bf Theorem 3.1.} {\em An abelian ring $R$ is VNL if and only if it
is an exchange ring such that for every idempotent $e$ of $R$,
either $eRe$ or $(1-e)R(1-e)$ is regular.}\\[3mm] {\bf
Proof.} The `only if' part follows from Lemma 2.4 and the
fact that every VNL ring is an exchange ring. Conversely, suppose
that $R$ is an abelian exchange ring such that for every
idempotent $e$ of $R$, either $eRe$ or $(1-e)R(1-e)$ is
regular. Let $a\in R$, then as $R$ is an exchange ring, there
exists an idempotent $e$ such that $e\in aR$ and $1-e\in (1-a)R$.
So $aR + (1-e)R = R$ and $eR + (1-a)R = R$ implying that  $eaR=eR$
and $(1-e)(1-a)R = (1-e)R$. Thus $ea$ and $(1-e)(1-a)$ are both
regular. Now if $eRe=eR$ is regular, then $e(1-a)$ is regular. So
$1-a=e(1-a)+(1-e)(1-a)$ is regular. Similarly, if $(1-e)R(1-e)$ is
regular, then $a$ is regular. $\Box$

\vspace{4mm} Recall that a ring $R$ is said to be potent if idempotents lift modulo $J(R)$ and every
right ideal not contained in $J(R)$ contains a non-zero idempotent. Every exchange ring is potent
 and a potent ring without infinite set of orthogonal idempotents is exchange. The following example shows
 that the above result is not true even for commutative potent rings.\\[3mm]
{\bf Example 3.2.} Let $R = \{(q_1,q_2,....,q_n,z,z,z, \dots):
n\geq 1 , q_i\in\mathbb{Q} ,z\in\mathbb{Z}\}$.  It is easy to see
that every non-zero ideal of $R$ contains a non-zero idempotent.
Also for any idempotent $e\in R$, either $eRe$ or
$(1-e)R(1-e)$ is regular but $R$ is not VNL as $\mathbb{Z}$, which
is a homomorphic image of $R$, is not VNL.

\vspace{3mm}
 In ([3], Lemma 2.7) Chen and Tong have shown that if $R$ is an abelian VNL ring, then $R/M(R)$ is a local ring.
 The following example shows that this may not be true for non-abelian rings:\\[3mm]
  {\bf Example 3.3.} Let $R=T_2(D)$, where $D$ is a division ring. By Corollary 2.13, $R$ is a VNL ring which
  is not regular. We will show that $M(R)=0$. Let
  $e= \left( \begin{array}{cc} a&b\\0&c \end{array} \right)$ be a non-zero idempotent in $R$. It is enough to show
  that $eR$ is not regular. It is clear that  $a, c \in \{0, 1\}$.
  If $a=c=1$, then $eR=R$. Also if $a=c=0$, then as $e$ is in $J(R)$, $b=0$. If $e= \left( \begin{array}{cc} 1&b\\0&0 \end{array} \right)$, then as $\left( \begin{array}{cc} 0&1\\0&0 \end{array} \right)$ is in $eR$, $eR$ is not a regular right ideal. Lastly if   $e= \left( \begin{array}{cc} 0&b\\0&1 \end{array} \right)$, then
   as $\left( \begin{array}{cc} 0&1\\0&0 \end{array} \right)$ is in $Re$, $Re$ is not a regular left ideal. Thus $M(R)=0$ and so $R/M(R)$ is not local.

   \vspace{4mm} The following lemma will give us another characterization of abelian VNL rings.\\[3mm]
 {\bf Lemma 3.4.} {\em Let $I$ be a regular ideal of a ring $R$. Then $R$ is VNL if and only if $R/I$ is VNL.}\\[3mm]
 {\bf Proof.} As any factor ring of a VNL ring is clearly VNL, we only have to prove the `if'part.
Suppose $R/I$ is VNL and  $a\in R$. Then either
$a+I$ or $1-a+I$ is regular
in $R/I$. In particular, either $a - axa\in I$ or
 $1-a -(1-a)y(1-a)\in I$ for some $x,y\in
R$. As $I$ is a regular ideal,  either $a-axa$  or
$(1-a)-(1-a)y(1-a)$ is a regular element of $R$. Thus by McCoy's
Lemma, either $a$ or $1-a$ is regular in $R$ showing that $R$ is
VNL.  $\Box$

\vspace{3mm} In view of ([3], Lemma 2.7) and above
Lemma, the following characterization of abelian VNL rings is
immediate.\\[3mm] {\bf Theorem 3.5.} {\em Let $R$ be an abelian
ring. Then $R$ is VNL if and only if $R/M(R)$ is a local
ring.}\\[7mm] {\bf 4. Characterization of semiperfect VNL rings}

\vspace{5mm}
 In this section we  characterize VNL rings without
infinite set of orthogonal idempotents, and also the VNL
 rings $R$ which have a primitive idempotent $e$ such that $eRe$
is not a division ring. As VNL rings without
infinite set of orthogonal idempotents are semiperfect,  we get a
characterization of semiperfect VNL rings.

Note that if $e$ is an idempotent in a ring $R$, then $$R\cong \left(\begin{array}{cc} S & X\\Y & T
\end{array}\right),$$ where $S=eRe$, $T=(1-e)R(1-e)$, $X=eR(1-e)$ is a $(S,T)$-bimodule and $Y=(1-e)Re$ is a $(T, S)$-bimodule such that $XY\subseteq S$ and $YX\subseteq T$. We will be tacitly using this representation of rings below and we will also be using $X$ and $Y$ in place of $XE_{12}$ and $YE_{21}$.\\[4mm]
{\bf Lemma 4.1.} {\em If $R= \left(\begin{array}{cc} S & X\\Y & T
\end{array}\right)$ such that $XY=0$ or $YX=0$, then $X,\;Y\subseteq J(R)$.}\\[3mm]
{\bf Proof.} If $XY=0$, then it is easy to see that $ \left(\begin{array}{cc} 0 & X\\0 & 0
\end{array}\right)$ is a quasi-regular right ideal of $R$ and $ \left(\begin{array}{cc} 0 & 0\\Y & 0
\end{array}\right)$ is a quasi-regular left ideal of $R$. Similarly if $YX=0$ then $\left(\begin{array}{cc} 0 & X\\0 & 0
\end{array}\right)$ is a quasi-regular left ideal of $R$ and $\left(\begin{array}{cc} 0 & 0\\Y & 0
\end{array}\right)$ is a quasi-regular right ideal of $R$. $\Box$
\\[4mm]{\bf Lemma 4.2.} {\em Let
$e_1$ and $e_2$ be two local idempotents of a ring. Then either
$e_1 R\cong e_2 R$ or $e_1 R e_2\subseteq J(R)$ and $
e_2Re_1\subseteq J(R)$}.\\[3mm] {\bf Proof}: Suppose $e_1 R\ncong e_2
R$. Then for any $r\in R$, $e_1 r e_2 R \neq e_1 R $. Because otherwise the map
from $e_2R\to e_1R$ given by the left multiplication with $e_1re_2$ splits implying that
$e_1 R\cong e_2 R$. Hence $e_1 r e_2 R$ is a
proper submodule of $e_1 R$, which has a unique maximal submodule
$e_1 J$. Thus $e_1 r e_2 R\subseteq e_1 J\subseteq J$ for every
$r$ implying that $e_1 R e_2\subseteq J(R)$. Similarly
$e_2Re_1\subseteq J(R)$. $\Box$\\[4mm]{\bf Corollary 4.3.} {\em A semiperfect
ring $R$ with $1 = e_1 + e_2$, where $e_1, e_2$ are orthogonal
primitive idempotents, is VNL if and only if $R$ is isomorphic to one of the
following:\\[2mm] (1) $M_2(D)$ for some division ring $D$.\\[2mm]
(2)  $\left(\begin{array}{cc} D & X\\ Y & L
\end{array}\right)$, where $D$ is a division ring, $L$ is a local
ring such that $XY =~0$.\\[2mm]
\indent In particular, if $J(R) = 0$,
then either $ R\cong M_2(D)$ or $ R\cong \left(\begin{array}{cc}
D_1 & 0\\ 0 & D_2 \end{array}\right),$ where $D_i$'s and $D$ are
division rings.}

\vspace{4mm}
\noindent {\bf Proof.} If $e_1 R\cong e_2 R$ then $R\cong
M_2(e_1 R e_1)$, where $e_1 R e_1$ is a local ring. So by
Lemma 2.4, $e_1 R e_1$ is a division ring. If $e_1 R\ncong e_2R$, then by Lemma 4.2,
$e_1 R e_2$ and $e_2 R e_1$ are contained in $J(R)$.  Again by Lemma 2.4,
either $e_1 R e_1$ or $e_2Re_2$ is a division ring.  We may assume that $e_1Re_1$ is a division ring.
Then $e_1 R e_2 R
e_1\subseteq e_1 R e_1\bigcap J(R) = e_1 J e_1 = 0$. So by taking $D=e_1Re_1$, $L=e_2Re_2$, $X=e_1Re_2$
 and $Y=e_2Re_1$, we see that $R\cong \left(\begin{array}{cc} D & X\\ Y & L
\end{array}\right)$ is as in (2) above.  $\Box$\\[4mm] {\bf Corollary 4.4.} {\em Let
$R$ be a semiperfect ring with $1 = e_1 + e_2 + e_3$, where
$e_1,e_2,e_3$ is an orthogonal set of primitive idempotents. Then
$R$ is VNL if and only if $R$ is one of the following:\\[2mm]
(1) $M_3(D)$ for some division ring $D$.\\[2mm]
(2) $\left(\begin{array}{cc} S & X\\ Y & L
\end{array}\right)$, where $S$ is semisimple and $L$ is a local  ring
with $XY = 0$.\\[2mm] (3) $\left(\begin{array}{cc} T & X\\ Y & D
\end{array}\right)$ where $T\cong$ $\left(\begin{array}{cc} D_1 &
X_1\\ Y_1 & D_2 \end{array}\right)$ is an NJ ring (see Example 2.1(2) above), $D$ a division
ring  with $YX =
0$}.\\[4mm] {\bf Proof.} If $e_iR\cong e_jR$ for all $i$, $j$,
 then $R\cong M_3(D)$ for some division ring $D$.

  It is clear that
$$R\cong \left(\begin{array}{cc}
(1-e_1)R(1-e_1) & (1-e_1)R e_1\\ e_1 R (1-e_1) & e_1 R e_1
\end{array}\right). \hspace{5mm} (A)$$

If $e_1 R e_1$ is local but not a division ring then $(1-e_1)R(1-e_1)$ is
a semisimple ring implying that $e_2 R e_2$ and $e_3 R e_3$ are division
rings. So by Lemma 4.2, $e_1 R e_2, e_2 R e_1, e_1 R e_3, e_3 R
e_1$ are all contained in $J(R)$ and so $(1-e_1)R e_1(1-e_1) = 0$. So $R$, in view  of (A),
 is as in (2) above.

 Now suppose all $e_i R e_i$ are division rings. If
 $e_2 R \cong e_3R$ but $e_1R \not \cong e_2R$, then $$(1-e_1)R(1-e_1)\cong M_2(D)$$ for some division ring
$D$ and, by Lemma 4.2, $$(1-e_1)R e_1 R (1-e_1) = 0 = e_1 R
(1-e_1)R e_1.$$ Thus $R$ as given in (A), is again as in (2) above.

 Lastly
assume that $e_1 R\ncong e_2 R\ncong e_3 R$. Then
$$(1-e_1)R(1-e_1)\cong \left(\begin{array}{cc} D_1 & X_1\\ Y_1 &
D_2 \end{array}\right),$$  with $X_1 Y_1 = 0 = Y_1 X_1$,
where $D_1=e_2Re_2$ and $D_3=e_3Re_3$.  In view of Lemma 4.2, it is clear that
 $ e_1 R (1-e_1)R e_1 = 0$. This in view of (A) implies that $R$
 is as in (3) above. $\Box$\\[4mm]
 {\bf Lemma 4.5.} {\em (1) Let $R=
 \left(\begin{array}{cc} S & X\\ Y & L \end{array}\right)$, where
$L$ is a local and  $S$ is a regular ring such that $XY = 0$. If $a=
\left(\begin{array}{cc} s & x\\y & u
\end{array}\right)\in R$ such that $u\in L$ is a unit, then $a$ is regular. In particular, $R$ is a VNL ring.\\[2mm]
(2) Let $S$ = $\left(\begin{array}{cc} T & X\\ Y & D
\end{array}\right)$, where $T\cong$ $\left(\begin{array}{cc} D_1 &
X_1\\ Y_1 & D_2 \end{array}\right)$ is an NJ ring, $D$ a division
ring with $YX = 0$. If $b=\left(\begin{array}{cc} t & x\\y & d
\end{array}\right)\in S$, then $b$ is regular under any of the following conditions:\\[2mm]
(a) If $t$ is regular in $T$ and $d\neq 0$.\\[2mm]
(b) If $t$ is a unit in $T$.}\\[3mm]
{\bf Proof.} Suppose $a=
\left(\begin{array}{cc} s & x\\y & u
\end{array}\right)\in R$ where $u\in L$ is a unit. Then clearly $$a = \left(\begin{array}{cc} s & xu^{-1}\\ 0 & 1 \end{array}\right) \left(\begin{array}{cc} 1&0\\y&u\end{array}\right)$$ As $\left(\begin{array}{cc} 1&0\\y&u\end{array}\right)$ is a unit in $R$ and $\left(\begin{array}{cc} s & xu^{-1}\\ 0 & 1 \end{array}\right)$ is regular by Proposition 2.8, $a$ is clearly regular.

Now suppose $$b=\left(\begin{array}{cc} t & x\\y & d
\end{array}\right)\in S,$$ where $t$ is regular in $T$ and $d\neq
0$. As all the regular elements of $T$ are unit regular, $t = eu$
for some idempotent $e$ and unit $u$ in $T$. Then $c =
\left(\begin{array}{cc} u & x\\y & d
\end{array}\right)$ is a unit in $S$ with $$c^{-1} = \left(\begin{array}{cc} (u-xd^{-1}y)^{-1} &
-u^{-1}xd^{-1}\\-d^{-1}y(u-xd^{-1}y)^{-1}& d^{-1}
\end{array}\right).$$ So $b$ is regular if and only if $bc^{-1}$
is regular. Now $$bc^{-1} = \left(\begin{array}{cc}
t_1 &
-exd^{-1}+xd^{-1}\\0 & 1
\end{array}\right),$$ where $t_1=t(u-xd^{-1}y)^{-1}-xd^{-1}y(u-xd^{-1}y)^{-1}$.  By Lemma 4.1, $x\in J(S)$ implying that $xd^{-1}y(u-xd^{-1}y)^{-1}\in
J(T)$. As $t(u-xd^{-1}y)^{-1}$ is regular in $T$ and only non-regular elements of $T$ are
the elements of $J(T)$, $t_1$ is regular. So by Proposition 2.8, $bc^{-1}$ is regular.

Lastly if $t$ is a unit and $d=0$, then again $b$ is regular with von
Neumann inverse as $\left(\begin{array}{cc} t^{-1} & 0\\0 & 0
\end{array}\right)$. $\Box$\\[4mm]\indent We are now ready to
characterize semiperfect VNL rings.\\[4mm]{\bf Theorem 4.6.} {\em
A semiperfect ring $R$ is VNL if and only if $R= A\times B$, where $A$ is a semisimple
ring and $B$ is one of the
following:\\[2mm] (1) Semisimple.\\[2mm](2) $R_1$ =
$\left(\begin{array}{cc} S & X\\ Y & L
\end{array}\right)$, where $L$ is a local ring, $S$ is a
semisimple ring such that $XY = 0$.\\[2mm] (3) $R_2$ =
$\left(\begin{array}{cc} T & X\\ Y & D
\end{array}\right)$, where $T\cong$ $\left(\begin{array}{cc} D_1 &
X_1\\ Y_1 & D_2 \end{array}\right)$ is an NJ ring, $D$ a division
ring such that $YX = 0$ (clearly this
case occurs in semiperfect rings with $1 = e_1+e_2+e_3$
only)}.\\[3mm] {\bf Proof.} In view of Lemma 4.5(1),  $R_1$ is VNL. Let $a=\left(\begin{array}{cc} t & x\\y & d
\end{array}\right)\in R_2$. If $t$ is regular in $T$ and $d\neq 0$, then $a$ is regular by Lemma 4.5(2). Also by Lemma 4.5(2), $a$ is regular if $t$ is a unit in $T$. Now assume that $t$ is not a unit in $T$. As non-regular elements
of $T$ are in $J(T)$,  $1-t$ is regular
 in $T$. If $d\neq 1$, then $1-a=\left(\begin{array}{cc} 1-t & -x\\-y & 1-d
\end{array}\right)$ which is regular by Lemma 4.5(2). Now suppose that $d=1$. Then if $t$ is regular in $T$, then
$a$ is regular. If $t$ is not regular, then $t\in J(T)$ and so $1-t$ is a unit in $T$. Then $1-a=\left(\begin{array}{cc} 1-t & -x\\-y & 0
\end{array}\right)$ is regular by Lemma 4.5 (2). Thus if $R\cong A\times B$, with $A$ semisimple and $B$ either
semisimple or isomorphic to $R_1$ or $R_2$, then $R$ is VNL.

 Conversely, let $R$ be a
semiperfect VNL ring. In view of the block decomposition of semiperfect rings and Lemma 2.4, $R\cong A \times B$, where $A$ is semisimple and $B$ is a semiperfect VNL ring with no non-trivial central idempotents. So we assume without loss of generality that $R$ is a semiperfect VNL ring without any non-trivial central idempotent.
In the proof below, we will call an
idempotent $e_i$, $1\leq i\leq n$, single if $e_i R\ncong e_j R$
for any $j\neq i$. Let $1=e_1+e_2+.....+e_n$ where $e_i$'s are
orthogonal primitive idempotents. We have already discussed the
case $n\leq 3$ in Corollary 4.3 and Corollary 4.4. We assume that
$n\geq 4$. Clearly for any $i$, $$R\cong \left(\begin{array}{cc}
(1-e_i)R(1-e_i)&(1-e_i)Re_i\\e_iR(1-e_i)&e_iRe_i\end{array}\right)
\hspace{8mm}(B) $$ Suppose first that there exist $e_i$ such that
$e_i R e_i$ is local but not a division ring. Then by Lemma 2.4,
$(1-e_i)R(1-e_i)$ is a semisimple ring. In particular, each $e_j R
e_j$ is a division ring whenever $j\neq i$. Thus $e_i R\ncong e_j
R$ and so by Lemma 4.2, $e_i R e_j, e_j R e_i$ are contained in
$J(R)$ for $j\neq i$. But as $(1-e_i)R(1-e_i)$ is semisimple,
$(1-e_i)R e_i R(1-e_i) = 0$. So $R$, in view of (B), is isomorphic
to $R_1$ in this case.

Now assume that each $e_i R e_i$ is a division ring. Note that if
$$\left(\begin{array}{cc} M_{n_1}(D_1) & X\\ Y & M_{n_2}(D_2)
\end{array}\right) \mbox { with } XY = 0 = YX$$ is VNL, then either it is
semisimple or one of $n_1$ or $n_2$ is equal to 1. So if each $e_i
R\cong e_j R$ for some $j\neq i$ then $R$ is a semisimple ring.
Now assume there exist $e_i$ such that $e_i R\ncong e_j R$ for any
$j$. If for each $e_j$, $i\neq j$, there exist $e_k$ with $k\neq j$ such that
$e_j R\cong e_k R$, then as mentioned above, $(1-e_i)R(1-e_i)$ is
a semisimple ring and by Lemma 4.2, $(1-e_i)Re_iR(1-e_i) = 0 =
e_iR(1-e_i)Re_i$. So $R$, in view of (B), is isomorphic to $R_1$
in this case also, with $L$ a division ring.\\ So we assume that
there are more than one single idempotents say $e_1, e_2,...
,e_r$. If $f_{r+1},\dots , f_m$ denote the sum of isomorphic
$e_i$'s. Then $$1 = e_1 +\dots + e_r+ f_{r+1}+\dots + f_m$$ Let
$e= f_{r+1}+\dots + f_m$ , by Lemma 2.4,  $eRe$ is clearly
semisimple. Suppose $(1-e)R(1-e)$ is also semisimple. If
$(e_i+f_j)R(e_i+f_j)$ is regular for every $i,j$,  then $R$ is
semisimple. We now assume that $(e_i+f_j)R(e_i+f_j)$ is not
regular for some $i, j$. If $f_j = e_{i_1} + e_{i_2} +\dots $ then
as $(e_i+e_{i_1})R(e_i+e_{i_1})$ is also not regular but
$(e_i+e_{i_1}+e_k+e_{i_2})R(e_i+e_{i_1}+e_k+e_{i_2})$ is VNL for each
$k\neq i$ (see Lemma 2.18), so by Lemma 2.4, $(e_k+e_{i_2})R(e_k+e_{i_2})$
is regular for each $k\neq i$ and hence $(e_k+f_j)R(e_k+f_j)$ is regular for each $k\neq i$. So $e_k R f_j = 0 = f_j R e_k$ for
each $k\neq i$. Thus $(1-e_i)R(1-e_i)$ is semisimple, this with
Lemma 4.2 implies $$(1-e_i)Re_i(1-e_i) = 0 = e_iR(1-e_i)Re_i$$ So
$R$, in view of (B), is isomorphic to $R_1$, with $L$ a division
ring.

 Lastly we assume that
$(1-e)R(1-e)$ is not semisimple. So there exist $e_i, e_j$ such
that $(e_i+e_j)R(e_i+e_j)$ is not semisimple and therefore by
Lemma 2.4, $(1-(e_i+e_j))R(1-(e_i+e_j))$ is semisimple. As
$n\geq 4$, we can pick  $k,l$ not equal to $i$ or $j$. Then by Lemma 2.18 and Lemma 2.4, either
$(e_k+e_i)R(e_k+e_i)$ or $(e_l+e_j)R(e_l+e_j)$ is semisimple.
Assume that $(e_k+e_i)R(e_k+e_i)$ is semisimple. If
$(e_k+e_j)R(e_k+e_j)$ is also semisimple and $e_k$ is single then
it is clearly central. If $e_k$ is not single, then the
corresponding $f_s$ is central. So assume that
$(e_k+e_j)R(e_k+e_j)$ is not semisimple. Then for any $t$ not in
\{i,j,k\}, $(e_t+e_i)R(e_t+e_i)$ is semisimple by Lemma 2.4
 and Lemma 2.18. In particular, $(1-e_j)R~(1-e_j)$ is semisimple. So by Lemma 4.2, $(1-e_j)Re_jR(1-e_j) = 0 = e_jR(1-e_j)Re_j$. So using $(B)$ with $i = j$, we have that $R$ is isomorphic to $R_1$,
  with $L$ a division ring. $\Box$\\[4mm]\indent A ring $R$ is called right
n-VNL-ring if $a_1 R + a_2 R +\dots + a_n R = R$ implies that some
$a_i$ is regular for some $i$. In
[4], it was shown that the semiperfect ring VNL ring $T_3(D)$ is
not 3-VNL. We prove that every semiperfect VNL ring is
$2$-VNL.\\[4mm] {\bf Theorem 4.7.} {\em A semiperfect VNL ring $R$ is
$2$-VNL.}\\[3mm]{\bf Proof.} We will use the characterization of semiperfect VNL rings as given in Theorem 4.6.
 We first show that $R_1$ is $2$-VNL. Suppose $A =
\left(\begin{array}{cc} s_1 & x_1\\y_1 & l_1
\end{array}\right)$ and $B = \left(\begin{array}{cc} s_2 &
x_2\\y_2 & l_2
\end{array}\right)$ are elements of $R_1$ such that $AR_1 + BR_1 = R_1$.
So there exist elements $C = \left(\begin{array}{cc} s_3 & x_3\\
 y_3 & l_3 \end{array}\right)$ and $D = \left(\begin{array}{cc} s_4 & x_4\\
 y_4 & l_4\end{array}\right)$ in $R_1$  such that $AC + BD = 1$ implying that
 $$\left(\begin{array}{cc} s_1s_3 + s_2s_4 & s_1x_3 + x_1l_3 + s_2x_4 + x_2l_4\\
 y_1s_3 + l_1y_3 + y_2s_4 + l_2y_4 & y_1x_3 + l_1l_3 + y_2x_4 + l_2l_4\end{array}\right) = \left(\begin{array}{cc} 1 & 0\\
 0 & 1\end{array}\right)$$ As $XY=0$ in $R_1$, by Lemma 4.1, $X, Y\subseteq J(R_1)$. Thus  $$\left(\begin{array}{cc} 0&0\\0&y_1x_3 + y_2x_4\end{array}\right) \in J(R_1)$$ and so $\left(\begin{array}{cc}1&0\\0&1-y_1x_3-y_2x_4\end{array}\right)$ is a unit in $R_1$ implying that $\left(\begin{array}{cc}1&0\\0&l_1l_3 + l_2l_4\end{array}\right)$ is a unit in $R_1$. So $l_1l_3 + l_2l_4$ is a unit in $L$. As $L$ is a local ring, either $l_1$ or $l_2$ is a unit in $L$. So in view of Lemma 4.5, either $A$ or $B$ is regular. Thus $R_1$ is $2$-VNL. We now show that $R_2$ is $2$-VNL. Suppose $P = \left(\begin{array}{cc} t_1&x_1\\y_1&d_1\end{array}\right)$ and $Q = \left(\begin{array}{cc} t_2&x_2\\y_2&d_2\end{array}\right)$ are elements of $R_2$ such that $PR_2 + QR_2 = R_2$. If $t_1$ or $t_2$ is a unit in $T$, then in view of Lemma 4.5, the corresponding element $P$ or $Q$ is regular in $R_2$. So suppose neither $t_1$ and nor $t_2$ is unit in $T$. As $PR_2 + QR_2 = R_2$, there exist $U = \left(\begin{array}{cc} t_3&x_3\\y_3&d_3\end{array}\right)$ and $V = \left(\begin{array}{cc} t_4&x_4\\y_4&d_4\end{array}\right)$ in $R_2$ such that $PU + QV = 1$, implying that $$\left(\begin{array}{cc} t_1t_3 + t_2t_4 + x_1y_3 + x_2y_4 & t_1x_3 + x_1d_3 + t_2x_4 + x_2d_4\\y_1t_3 + d_1y_3 + y_2t_4 + d_2y_4 & d_1d_3 + d_2d_4\end{array}\right) = \left(\begin{array}{cc} 1&0\\0&1\end{array}\right)$$ So $d_1d_3 + d_2d_4 = 1$, implying that at least one of $d_1$ and $d_2$ is a unit. Again as $YX=0$ in $R_2$, by Lemma 4.1, $X, Y\subseteq J(R_2)$ and so $\left(\begin{array}{cc} x_1y_3 + x_2y_4&0\\0&0\end{array}\right)\in J(R_2)$, $\left(\begin{array}{cc} 1-x_1y_3-x_2y_4&0\\0&1\end{array}\right)$ is a unit in $R_2$. So $t_1t_3 + t_2t_4$ is a unit in $T$. As $T$ is an NJ ring and none of $t_1$ and $t_2$ is a unit, it is easy to see that both $t_1$ and $t_2$ are regular in $T$. Also since one of $d_1$ and $d_2$ is a unit in $D$, so by Lemma 4.5, the corresponding element $P$ or $Q$ is regular in $R_2$. $\Box$\\[4mm]        \indent We now characterize VNL rings in which there
is a primitive idempotent $e$ such that $eRe$ is not a division
ring.\\[4mm] {\bf Theorem 4.8.} {\em Let $R$ be a ring
with a primitive idempotent $e$ such that $eRe$ is not a division
ring. Then $R$ is VNL if and only if $R\cong$
$\left(\begin{array}{cc} S & X\\ Y & L \end{array}\right)$, where
$L$ is a local ring, $S$ is a regular ring and $XY =
0$.}\\[3mm]
{\bf Proof.} The `if' part follows from Lemma 4.5(1). Now suppose that $R$ is VNL and $e\in R$ is a primitive idempotent such that $eRe$ is not a division ring. So $eRe$ is a local ring and, by Lemma 2.4, $(1-e)R(1-e)$ is a regular ring.  We have $$R\cong\left(\begin{array}{cc} (1-e)R(1-e) & (1-e)Re\\
eR(1-e) & eRe
\end{array}\right).$$  We now show that $eR(1-e) \subseteq J(R)$. Note that if for any element $r\in R$, $er(1-e)R=eR$, then $eR$ will be isomorphic to a summand of $(1-e)R$. But as corner rings of regular rings are regular, $eRe$ is regular and hence a division ring, a contradiction. So $er(1-e)R$ is a proper submodule of a local module $eR$ implying that $er(1-e)R \subseteq eJ(R)$. So $eR(1-e) \subseteq J(R)$. In particular, $(1-e)ReR(1-e)\in J(R)\cap (1-e)R(1-e)=0$, as $(1-e)R(1-e)$ is a
regular ring.  $\Box$\\[8mm]\noindent {\bf 5. A Sufficient Condition}

 \vspace{5mm}
  A ring $R$ is called semipotent if every right ideal not contained in $J(R)$ contains a nonzero idempotent. In general we have
  $$ \mbox{ NJ } \Longrightarrow \mbox{ VNL } \Longrightarrow \mbox{ Exchange } \Longrightarrow \mbox{ Potent  } \Longrightarrow \mbox{ Semipotent},$$
  with none of the implications reversible. We give below a sufficient condition for all these classes of rings to
  coincide.\\[4mm] {\bf Lemma 5.1.} {\em Let $R$ be a semipotent ring without central
idempotents such that $J(R)\neq 0$ but $J(eRe)=0$ for every proper
idempotent $e$ of $R$. Then the following hold:\\[2mm](1) For every
proper idempotent  $e$ of $R$, $eR(1-e)$ and $(1-e)Re$ are
contained in $J(R)$.\\[2mm](2) If $0\neq e = e^2$ is such that $ae = a
= ea$ for every $a$ in $J(R)$, then $e$ is central and hence
$e=1$}
\\[3mm]{\bf Proof.} First note that if $a$ is in $J(R)$, then
for every proper idempotent $e$ of $R$, $a = ea + ae$ as
$(1-e)J(1-e)= 0$ and $eJe = 0$. Let $e$ be proper idempotent of
$R$. Now $(1-e)JeR(1-e) = 0$ as it is contained in $(1-e)J(1-e)$,
so $JeR(1-e) = eJeR(1-e) = 0$. Also $eR(1-e)J(1-e) = 0$ implying
that $eR(1-e)J = eR(1-e)Je = 0$. Thus $eR(1-e)J = JeR(1-e) = 0$
and so $eR(1-e)\subseteq Ann(J(R))$. If $eR(1-e)\nsubseteq J(R)$
then Ann$J(R)$ $\nsubseteq$ $J(R)$ and as $R$ is semipotent, there
exist $0\neq f=f^2\in AnnJ(R)$. Then as $J(R)\neq 0$, $f\neq 1$
and hence $f$ is proper. So $a = af+fa = 0$ for every $a$ in
$J(R)$ implying $J(R) = 0$ , a contradiction. Hence
$eR(1-e)\subseteq J(R)$. Similarly $(1-e)Re\subseteq J(R)$.

 Now
suppose $0\neq e=e^2$ is such that $ae = ea = a$ for every $a$ in
$J(R)$. If $e=1$, then nothing to prove. If $e$ is proper, then as
$eR(1-e)$ and $(1-e)Re$ are contained in $J(R)$, $e.er(1-e) =
er(1-e)e = 0$ for every $r$ in $R$ implying that $er(1-e) = 0$,
similarly $(1-e)re = 0$ for all $r$ in  $R$ and so $e$ is central.~$\Box$\\[4mm]{\bf Proposition 5.2.} {\em Let $R$ be a ring without
central idempotents such that $J(R)\neq 0$ but $J(eRe) = 0$ for
every proper idempotent $e$ of $R$. Then the following are
equivalent:\\[2mm] (1) $R$ is a VNL ring.\\ (2) $R$ is an exchange
ring.\\ (3) $R$ is a potent ring.\\ (4) $R$ is a semipotent
ring.\\ (5) $R$ is an NJ ring.} \\[3mm] {\bf Proof.} The
implications (1) $\Longrightarrow$ (2) $\Longrightarrow$ (3)
$\Longrightarrow$ (4) and (5) $\Longrightarrow$ (1) hold in
general. So we only have to prove the implication (4) $\Longrightarrow$ (5).

 Note
that if $a$ is in $J(R)$, then for every proper idempotent $e$ of
$R$, $a = ea + ae$. Now we prove that if $e$ is a proper
idempotent of $R$ then $eRe$ has only trivial idempotents. Suppose
$f$ is a proper idempotent in $eRe$, then $ef = fe = f$. Clearly
$e-f\neq 0$, also it is easy to see that $e-f\neq 1$ and hence
$e-f$ is a proper idempotent of $R$. Now for any $a\in J(R)$, $a =
af+fa$. Also $af = af(e-f) + (e-f)af = eaf$ and $fa = fa(e-f)+
(e-f)fa = fae$. Then as $a = af+fa$, $ea = eaf+efa = af+fa = a$
and $ae = afe+fae = af+fa = a$. So we have $ea = ae = a$ for every
$a$ in $J(R)$ and thus by Lemma 5.1, $e=1$, a contradiction. Thus
$eRe$ has only two idempotents. Now as $R$ semipotent implies
$eRe$ semipotent, $eRe$ is a local ring for every proper
idempotent $e$ of $R$. Also as $J(eRe) = 0 $, $eRe$ is a division
ring.

 Now if $R$ has no proper idempotent then $R$ is local and
hence NJ. If $e$ is a proper idempotent in $R$, then $eRe$ and
$(1-e)R(1-e)$ are division rings and in view of Lemma 5.1,
$eR(1-e)Re = 0 = (1-e)ReR(1-e)$. Thus $$R\cong
\left(\begin{array}{cc} D_1 & X\\ Y & D_2
\end{array}\right),$$ where $D_1, D_2$ are division rings and $XY =
0 = YX$ and hence by Nicholson's characterization of NJ rings,
$R$ is an NJ ring (see [6]). $\Box$\\[3mm] \indent In Theorem
3.1, we proved that an abelian ring $R$ is VNL if and only if it
is an exchange ring with the property that for every idempotent
$e$ of $R$, one of the two corner rings $eRe$ or $(1-e)R(1-e)$ is
regular. So one may ask:\\[3mm] {\bf Question 5.3.}  {\em Let $R$ be an arbitrary exchange ring
with the property that for every idempotent $e$ of $R$, one
of $eRe$ and $(1-e)R(1-e)$ is regular. Then is $R$ a VNL ring?}\\[3mm] \indent
In Theorem 4.7, we proved that a semiperfect VNL ring is $2$-VNL.
So the following natural question arises:\\[3mm] {\bf Question 5.4} {\em Is
every VNL ring $2$-VNL?}
\\[8mm]{\bf References}

\vspace{3mm}\noindent [1] \noindent B. Brown and N. McCoy, {\em
The Maximal Regular Ideal of a Ring}, Proc. Amer. Math.
Soc. {\bf 1}(1950), 165-171.\\[2mm] [2] \noindent V. P. Camillo, H.
P. Yu , {\em Exchange Rings, Units and Idempotents},
Comm. Algebra {\bf 22}(1994), 4737-4749.\\[2mm] [3] \noindent W.
X. Chen and W. T. Tong, {\em On Noncommutative VNL-rings and GVNL-rings},
Glasgow Math. J., {\bf 48}(2006),
11-17.\\[2mm] [4] \noindent J. Chen and Z. Ying, {\em On VNL-rings
and n-VNL-rings}, to appear.\\[2mm] [5] \noindent  M. Contessa,
{\em On Certain Classes Of PM-rings}, Communications in
Algebra, {\bf 12}(1984), 1447-1469.\\[2mm] [6] \noindent W. K.
Nicholson, {\em Rings whose elements are Quasi-Regular or
Regular}, Aequationes Mathematicae {\bf 9}(1973),
64-70.\\[2mm] [7] \noindent W. K. Nicholson, {\em Lifting
Idempotents and Exchange Rings}, Trans. Amer. Math. Soc.
{\bf 229}(1977), 269-278.\\[2mm] [8] \noindent E. A. Osba, M.
Henriksen and O. Alkam, {\em Combining Local and von
Neumann regular rings}, Communications in Algebra, {\bf
32}(2004), 2639-2653.\\[2mm] [9] \noindent  E. A. Osba,
M. Henriksen, O. Alkam and F. A. Smith, {\em The
Maximal Regular Ideal of Some Commutative Rings}, Comment. Math.
Univ. Carolinae, {\bf 47}(2006), 1-10.

\end{document}